\documentclass[a4paper,12pt,twoside]{article}
\usepackage{amsmath,amssymb,amsthm,verbatim}

\numberwithin{equation}{section}
\newtheorem{thm}[equation]{Theorem}

\newtheorem{lem}[equation]{Lemma}

\newcounter{mycount}
\newenvironment{romlist}{\begin{list}{\rm(\roman{mycount})}%
   {\usecounter{mycount}\labelwidth=1cm\itemsep 0pt}}{\end{list}}
\newenvironment{numlist}{\begin{list}{(\arabic{mycount})}%
   {\usecounter{mycount}\labelwidth=1cm\itemsep 0pt}}{\end{list}}

{\noindent{\sc Proof. }}{{\hfill $\Box$}\par\vskip2\parsep}
{\par\vskip2\parsep\noindent{\sc Proof of\ #1. }}{{\hfill $\Box$}
  \par\vskip2\parsep}
{\noindent{\sc Sketch Proof. }}{{\hfill $\Box$}\par\vskip2\parsep}
{\par\vskip2\parsep\noindent{\sc Sketch Proof of\ #1. }}{{\hfill $\Box$}\par
  \vskip2\parsep}

\def\o{{\mathrm o}}

\def\s{\sigma}
\def\qq{\qquad}
\def\q{\quad}

\def\rc{random-cluster}

\def\ZZ{{\mathbb Z}}
\def\RR{{\mathbb R}}
\def\EE{{\mathbb E}}
\def\LL{{\mathbb L}}
\def\PP{{\mathbb P}}

\def\la{\langle}

\def\Om{\Omega}
\def\om{\omega}

\def\oo{\infty}

\def\rad{\text{rad}}

\def\sF{\mathcal F}

\def\a{\alpha}

\def\f{\rho}

\def\psd{{p_{\mathrm{sd}}}}
\def\pc{p_{\mathrm c}}
\def\coin{\text{conditional influence}}
\def\ain{\text{absolute influence}}
\def\be{\begin{equation}}
\def\ee{\end{equation}}
\def\eps{\epsilon}
\def\bigmid{\,\big|\,}
\def\cov{\mathrm{cov}}
\def\q{\quad}
\def\qq{\qquad}
\def\comp#1{\overline {#1}}
\def\lest{\le_{\mathrm{st}}}
\def\sm{\setminus}
\def\fpq{\phi_{p,q}}
\def\resp{respectively}
\def\bx{\mathbf x}
\def\by{\mathbf y}

\def\la{\lambda}
\def\ol#1{\overline{#1}}

\def\LW{\mathrm{LW}}
\def\sA{\mathcal A}
\def\sH{\mathcal H}
\def\sV{\mathcal V}
\def\SW{\mathrm{SW}}

\title{Influence and sharp-threshold theorems for monotonic measures}
\author{B.\ T.\ Graham, G.\ R.\ Grimmett\\
   {\small Statistical Laboratory, University of Cambridge,}\\
   {\small Wilberforce Road, Cambridge CB3 0WB, U.K.}\\
   {\small 3 May 2005}}

\begin{document}
\maketitle

\begin{abstract}The influence theorem for product
measures on the discrete space $\{0,1\}^N$
may be extended to probability measures with
the property of monotonicity (which is equivalent to
`strong positive-association'). Corresponding
results are valid for probability measures
on the cube $[0,1]^N$ that are absolutely
continuous with respect to Lebesgue measure.
These results lead to
a sharp-threshold theorem for measures of random-cluster type, 
and this may be applied to box-crossings 
in the two-dimensional \rc\ model.
\renewcommand{\thefootnote}{}
\renewcommand{\thefootnote}{\arabic{footnote}}
\setcounter{footnote}{0}
\end{abstract}

\section{Introduction}
Influence and sharp-threshold theorems have proved useful in
the study of 
problems in discrete probability. Reliability theory and random graphs
provided early problems of this type, followed by percolation.
Important progress has been made since \cite{BOL, KKL} towards
a general theory, of which one striking aspect has been the
use of discrete Fourier analysis and hypercontractivity.
The reader is referred to \cite{F04, FKST} for a history
and bibliography.

Let $\Om=\{0,1\}^N$ where $N<\oo$,
and let $\mu_p$ be the product measure on $\Om$ with
density $p$. Vectors in $\Om$ are denoted
by $\om=(\om(i): 1\le i\le N)$.
For any increasing 
subset $A$ of $\Om$, and any $i\in\{1,2,\dots,N\}$,
we define the {\it\coin\/} $I_A(i)$ by
\be
I_A(i) = \mu_p(A\mid X_i=1) - \mu_p(A\mid X_i=0),
\label{inf}
\ee
where $X_i$ is the indicator function of the event
$\{\om\in\Om: \om(i)=1\}$.
It is well known (see \cite{BKKKL, FKST, KKL, Tal94}) that there 
exists an absolute positive constant $c$ such that the following
holds. For
all $N$, all $p\in(0,1)$,
and all increasing $A$, there exists $i\in \{1,2,\dots,N\}$
such that 
\be
I_A(i) \ge c\min\{\mu_p(A), 1-\mu_p(A)\} \frac{\log N}N.
\label{inf0}
\ee
The proof uses discrete Fourier analysis
and a technique known as `hypercontractivity'.
Inequality (\ref{inf0}) is usually stated for
the case $p=\frac12$, but it holds with the same 
constant $c$ for
all $p\in(0,1)$.

There is an important application to the theory of sharp
thresholds for product measures.  Let $\Pi_N$ be the set of all 
permutations of the index set 
$I=\{1,2,\dots,N\}$. A subgroup $\sA$ of $\Pi_N$
is said to {\it act transitively\/} on $I$ if, for
all distinct
pairs
$j,k\in I$, there exists $\pi\in\sA$
with $\pi_j=k$.  Any $\pi\in\Pi_N$ acts on $\Om$ by
$\pi\om=(\om(\pi_i): 1\le i\le N)$.
An event $A$ is called {\it symmetric\/}
if there exists a subgroup $\sA$ of $\Pi_N$ acting transitively
on $I$ such that $A=\pi A$.
If $A$ is symmetric, then $I_A(j)=I_A(k)$
for all $j$, $k$. By summing (\ref{inf0}) over $i$ we obtain
for symmetric $A$ that
\be
\sum_{i=1}^N I_A(i) \ge c\min\{\mu_p(A), 1-\mu_p(A)\} \log N.
\label{suminf}
\ee
It is standard (see the discussion of Russo's formula
in \cite{G99}) that
\be
\frac{d}{dp} \mu_p(A) = \sum_{i=1}^N I_A(i),
\label{Russo}
\ee
and it follows as in \cite{FKST} that, for $0<\eps<\frac12$,
the function $f(p)=\mu_p(A)$ 
increases from $\eps$ to $1-\eps$ over an interval
of values of $p$ with length smaller in order
than $1/\log N$.

We refer to such a statement as a `sharp-threshold theorem',
and we note that such results have wide applications
to problems of discrete probability.
For example, the observations above
have been used recently 
in \cite{BR1} to obtain a further
proof of the
famous theorem of Harris and Kesten that the critical probability
$\pc$ of bond percolation on the square
lattice satisfies $\pc=\frac12$. 
Using a similar argument in a second paper,
\cite{BR2},
they have proved the conjecture that the critical probability
of site percolation on a certain Poisson--Voronoi graph in $\RR^2$
equals $\frac12$
almost surely.

The principal purpose of the current article is to extend the results
above to probability measures more general than
product measures. We shall prove such results
for measures  having a certain condition of
`monotonicity', which
is equivalent to the FKG lattice condition
and is described in the next section.
There are many situations in the probabilistic
theory of statistical mechanics where such measures
are encountered, including the Ising model and the
random-cluster model.

We define monotonic probability measures in Section 2,
and  we note there that monotonicity is equivalent to the FKG lattice condition.
This is followed by an influence theorem for monotonic measures. 

A monotonic measure $\mu$ may be used as the basis of a certain parametric
family of measures on $\Om$ indexed by a parameter $p\in(0,1)$.
The influence theorem for $\mu$  may then be used to obtain 
a sharp-threshold theorem for this class, as described in Section 3.

The influence theorem on the discrete space
$\Om$ was extended in \cite{BKKKL} to
product measures on the Euclidean cube $[0,1]^N$.
Using the methods of Section 2, similar
results may be proved for general monotonic measures on $[0,1]^N$.
Unlike the discrete case, such an influence theorem does not 
appear to imply a corresponding sharp-threshold theorem.
This is discussed in Section 4.

We turn finally to the \rc\ model, which may be viewed as an
extension of percolation and a generalization of the Ising/Potts models for
ferromagnetism, see \cite{G03, G-RC}.  
The \rc\ measure is defined in Section 5,
and the sharp-threshold theorem is applied to the existence
of box-crossings in two dimensions.

\section{Influence for monotonic measures}
We begin this section with a classification,
further details of which may be found in \cite{G-RC}.
Let $1\le N<\oo$, and write $I=\{1,2,\dots,N\}$
and $\Om=\{0,1\}^N$. The set of all subsets of
$\Om$ is denoted by $\sF$.
A probability measure $\mu$ on $(\Om,\sF)$ is said to
be {\it positive\/} if $\mu(\om)>0$ for
all $\om\in\Om$. It is said to satisfy the
{\it FKG lattice condition\/} if
\be
\mu(\om_1\vee \om_2)\mu(\om_1\wedge\om_2)\ge \mu(\om_1)\mu(\om_2)
\qq\text{for all } \om_1,\om_2\in\Om,
\label{FKG}
\ee
where
$\om_1\vee\om_2$ and
$\om_1\wedge\om_2$ are given by
\begin{align*}
\om_1\vee\om_2(i) &=\max\{\om_1(i),\om_2(i)\},   &i\in I, \\
\om_1\wedge\om_2(i) &=\min\{\om_1(i),\om_2(i)\}, &i\in I.
\end{align*}
See \cite{FKG, G-RC}.

The set $\Om$ is a partially ordered set with
the partial order:
$\om\ge\om'$ if $\om(i)\ge \om'(i)$
for all $i\in I$.
A non-empty event $A \in \sF$ is called {\it increasing\/} if:
$\om\in A$ whenever there exists $\om'$ with 
$\om\ge\om'$ and $\om'\in A$. It is called
{\it decreasing\/} if its complement is increasing.
For probability measures $\mu_1$, $\mu_2$ on
$(\Om,\sF)$, we write $\mu_1 \lest \mu_2$, and say that
$\mu_1$ is dominated stochastically by $\mu_2$, if
$$
\mu_1(A) \le \mu_2(A)\qq\text{for all increasing events } A.
$$
The indicator function of an event $A$ is denoted by $1_A$.
For $i\in I$, we write $X_i$ for the indicator function
of the event $\{\om\in\Om: \om(i)=1\}$.

A probability measure $\mu$ on $\Om$ 
is said to be {\it positively associated\/}
if
$$
\mu(A\cap B) \ge \mu(A)\mu(B)\qq
\text{for all increasing events $A$, $B$}.
$$
The famous FKG inequality of \cite{FKG}
asserts that a positive probability measure
$\mu$ is positively associated if it satisfies 
the FKG lattice condition. It is well known that the FKG lattice condition
is not necessary for positive association, and we explore this next.

We shall for simplicity restrict ourselves henceforth to
positive measures.
The FKG lattice condition is equivalent to a 
stronger property termed `strong positive association'. 
For $J\subseteq I$ and $\xi\in\Om$, let
$\Om_J=\{0,1\}^J$ and 
\be
\Om_J^\xi=\{\om\in\Om:\om(j)=\xi(j) \text{ for } j \in I\sm J\}.
\label{new2.2}
\ee
The set of all subsets of $\Om_J$ is denoted by $\sF_J$.
Let $\mu$ be a positive probability measure on $(\Om,\sF)$, and define the
conditional probability measure $\mu_J^\xi$ on $(\Om_J,\sF_J)$  by
\be
\mu_J^\xi(\om_J) = 
\mu\bigl(X_j=\om_J(j) \text{ for } j \in J\bigmid 
X_i=\xi(i) \text{ for }
i\in I\sm J\bigr),
\qq \om_J\in\Om_J.
\label{new2.3}
\ee
We say that $\mu$ is {\it strongly positively-associated\/}
if: 
for all $J\subseteq I$ and all $\xi\in\Om$, the measure
$\mu_J^\xi$ is positively associated.

We call $\mu$ {\it monotonic\/} if: for all $J\subseteq I$, all
increasing subsets $A$ of $\Om_J$, and
all $\xi,\zeta\in\Om$,
\be
\mu_J^\xi(A) \le 
\mu_J^\zeta(A)\qq\text{whenever } \xi\le\zeta.
\label{new2.3a}
\ee
That is, $\mu$ is monotonic if, for all $J\subseteq I$, 
\be
\mu_J^\xi\lest\mu_J^\zeta \qq\text{whenever } \xi\le\zeta.
\label{new2.4}
\ee
We call $\mu$ 1-{\it monotonic\/} if (\ref{new2.4}) holds for all singleton
sets $J$. That is, $\mu$ is 1-monotonic if and only if, for all $j\in I$,
\be
\mu\bigl(X_j=1 \bigmid X_i=\xi(i) \text{ for all } i \in I\sm\{j\}\bigr)
\label{new2.5}
\ee
is non-decreasing in $\xi$.

The following theorem is fairly standard, and the proof may be
found in \cite{G-RC}.

\begin{thm}\label{new2.7}
Let $\mu$ be a positive probability measure on $(\Om,\sF)$.
The following are equivalent.
\begin{romlist}
\item $\mu$
is strongly positively-associated.
\item $\mu$ 
satisfies the FKG lattice condition.
\item $\mu$ is monotonic.
\item $\mu$ is $1$-monotonic.
\end{romlist}
\end{thm}

Our principal influence theorem is as follows. For a positive
probability measure $\mu$ and an increasing
event $A$, the 
{\it\coin\/} of the index $i$ ($\in I$)
is given as in (\ref{inf}) by
\be
I_A(i) = \mu(A\mid X_i=1) - \mu(A\mid X_i=0).
\label{inf2}
\ee
For a product measure $\mu_p$, the influence of 
the index $i$ was defined in \cite{BOL, KKL} 
as $\mu_p(\om^i\in A,\ \om_i\notin A)$, where $\om^i$
(\resp, $\om_i$) denotes the configuration obtained from $\om$ by
setting $\om(i)$ equal to 1 (\resp, 0).
We refer to the latter quantity as the {\it\ain\/} of index $i$.
The absolute and conditional influences are equal for product measures,
but one should note that 
\be
I_A(i) \ne \mu(\om^i\in A,\ \om_i\notin A)
\label{inf3}
\ee
for general probability measures $\mu$. 
Further discussion of this point is provided after the next theorem.

\begin{thm}[Influence]\label{infthm}
There exists a constant $c$ satisfying $c\in(0,\oo)$
such that the following holds. Let $N\ge 1$ and let $A$
be an increasing
subset of $\Om=\{0,1\}^N$. 
Let $\mu$ be a positive probability
measure on $(\Om,\sF)$ that is monotonic. 
There exists $i\in I$ such that
\be
I_A(i) \ge c\min\{\mu(A),1-\mu(A)\} \frac{\log N}N.
\ee
\end{thm}

Since product measures are monotonic, this extends the influence
theorem of \cite{KKL}. In the proof of Theorem \ref{infthm},
we shall encode the measure $\mu$ in terms of Lebesgue measure on
$[0,1]^N$, and we shall appeal to the influence theorem of \cite{BKKKL}.
Thus, we shall require no further arguments of discrete
Fourier analysis than those already present in \cite{BKKKL, KKL}.

We return briefly to the discussion of
absolute and conditional influences. Suppose, for
illustration, that $P$ is chosen at random with
$\PP(P=\frac13)=\PP(P=\frac23)=\frac12$ and that,
conditional on the value of $P$, we are provided with
independent Bernoulli random
variables $X_1,X_2,\dots,X_N$ with parameter $P$. Consider the increasing
event $A=\{S_N > \frac12 N\}$, where $S_N=X_1+X_2+\dots+X_N$.
By symmetry, the \coin\ of each index 
is the same, as is the \ain\ of each index.
It is an easy calculation that
$$
I_A(1) = \tfrac13 + \o(1)\qq\text{as } N\to\oo.
$$
On the other hand,
\begin{align*}
\PP(\om^1\in A,\ \om_1\notin A)&=
\PP\left(\tfrac12 N - 1 < \sum_{i=2}^N X_i \le \tfrac12 N\right)\\
&=\o(e^{-\gamma N})\qq\text{as }N\to\oo,
\end{align*}
for some $\gamma>0$. This example indicates not only
that the absolute and conditional influences can
be very different, but
also that the conclusion of Theorem \ref{infthm} would be false if
re-stated for \ain s.

In the proof of Theorem \ref{infthm} following,
we see that monotonicity has the effect of increasing
the influence of each coordinate in $I$.

\begin{proof}[Proof of Theorem \ref{infthm}]
Let $A \in\sF$ be an increasing event, and let $\mu$ be
positive and monotonic. Let $\la$ denote 
Lebesgue measure on the cube $[0,1]^N$.
We propose to construct an increasing
subset $B$ of $[0,1]^N$ with the property that $\la(B)=\mu(A)$,
to apply the influence theorem of \cite{BKKKL} to the set $B$, and to deduce
the claim. This will be done via a certain function
$f:[0,1]^N \to \{0,1\}^N$ that we construct next.

Let $\bx=(x_i: 1\le i \le N)\in [0,1]^N$, and let 
$f(\bx)=(f_i(\bx): 1\le i\le N)$ be
given recursively as follows. The first coordinate
$f_1(\bx)$ is defined by:
\be
\text{with}\q a_1=\mu(X_1=1),\q\text{set} \q f_1(\bx) =
\begin{cases}
 1 &\text{if } x_1 > 1-a_1,\\
0 &\text{otherwise}.
\end{cases}
\label{a1}
\ee
Suppose we know $f_i(\bx)$ for $1\le i < k$. Let
\be
a_k=\mu(X_k=1\mid X_i = f_i(\bx)\text{ for } 1\le i< k),
\label{a2}
\ee
and define 
\be
f_k(\bx) =\begin{cases}
 1 &\text{if } x_k > 1-a_k,\\
0 &\text{otherwise}.
\end{cases}
\label{a2'}
\ee

Suppose that $\bx\le \bx'$, and write $a_k$ and $a_k'$ for
the corresponding values in (\ref{a1})--(\ref{a2}). Clearly
$a_1=a_1'$, so that $f_1(\bx) \le f_1(\bx')$. Since $\mu$
is monotonic, $a_2\le a_2'$, so that $f_2(\bx)\le f_2(\bx')$.
Continuing inductively, we find that $f_k(\bx)\le f_k(\bx')$
for all $k$, which is to say that $f(\bx)\le f(\bx')$. Therefore,
$f$ is non-decreasing on $[0,1]^N$.
Let $B$ be the increasing subset of $[0,1]^N$ given by
$B=f^{-1}(A)$. 

We make four notes concerning the definition of $f$.
\begin{numlist}
\item Each $a_k$ depends only on
$x_1,x_2,\dots,x_{k-1}$.
\item Since $\mu$ is positive, the $a_k$
satisfy $0<a_k<1$ for all $\bx\in[0,1]^N$ and $k\in I$.
\item For any $\bx\in [0,1]^N$ and $k\in I$,
the values $f_k(\bx),f_{k+1}(\bx),\dots,f_N(\bx)$
depend on $x_1,x_2,\dots,x_{k-1}$ only through
the values $f_1(\bx),f_2(\bx),\dots,f_{k-1}(\bx)$.
\item The function $f$ and the event $B$ depend on
the ordering of the set $I$.
\end{numlist}

Let $U=(U_i: 1\le i\le N)$ be the identity function on
$[0,1]^N$, and note that $U$ has law $\la$. By the method
of construction
of the function $f$, $f(U)$ has law $\mu$.
In particular,
\be
\mu(A) = \la(f(U)\in A) = \la(U\in f^{-1}(A)) = \la(B).
\label{a10}
\ee

Let
$$
J_B(i) = \la(B\mid U_i=1)- \la(B\mid U_i=0),
$$
where the conditional probabilities are to be interpreted as
$$
\la(B\mid U_i=u) = \lim_{\eps\downarrow 0}\left\{
\frac 1\eps \la(B\mid U_i\in(u-\eps,u+\eps))\right\},\qq u=0,1.
$$
Since $B$ is an event with a certain simple structure, this
is the same as $\la_{N-1}(B_i^u)$ for $u=0,1$, 
where $\la_{N-1}$ is $(N-1)$-dimensional
Lebesgue measure and 
$B_i^u$ is the set of all $(N-1)$-vectors
$(x_1,\dots,x_{i-1},x_{i+1},\dots,x_N)$ such that
$(x_1,\dots,x_{i-1},u,x_{i+1},\dots,x_N)\in B$.

By Theorem 1 of \cite{BKKKL}, we may find a constant $c>0$,
independent of
the choice of $N$ and $A$, such that: there exists $i\in I$ with 
\be
J_B(i) \ge c \min\{\la(B),1-\la(B)\}\frac{\log N}N.
\label{a4}
\ee
We choose $i$ accordingly.

We claim that
\be
I_A(j) \ge J_B(j)\qq\text{for } j\in I.
\label{a5}
\ee
Once (\ref{a5}) is shown, the claim follows from (\ref{a10})
and (\ref{a4}).
We prove next that
\be
I_A(1) \ge J_B(1).
\label{a6}
\ee
We have that
\begin{align}
I_A(1) &= \mu(A\mid X_1=1) - \mu(A\mid X_1=0)\nonumber\\
&= \la(B\mid f_1(U)=1) - \la(B\mid f_1(U)=0)\nonumber\\
&= \la(B\mid U_1 > 1-a_1) - \la(B\mid U_1\le 1-a_1)\nonumber\\
&= \la(B\mid U_1=1) - \la(B\mid U_1=0)\nonumber\\
&= J_B(1),
\label{a7}
\end{align}
where we have used notes (2) and (3) above.
This implies (\ref{a6}).

We turn our attention to (\ref{a5}) with $j\ge 2$. 
We re-order
the set $I$ to bring the index $j$ to the front. That is, we let
$K$ be the re-ordered index set $K=
(k_1,k_2,\dots,k_N) = (j,1,2,\dots,j-1,j+1,\dots,N)$.
We write $g=(g_{k_i}: 1\le i\le N)$ for the associated function given by
(\ref{a1})--(\ref{a2'}) subject to the new ordering,
and $C=g^{-1}(A)$.
Thinking of (\ref{a1})--(\ref{a2'}) as an algorithm
for constructing $f$, we are applying the same algorithm
to the re-ordered set $K$.

We claim that
\be
J_{C}(k_1) \ge J_B(j).
\label{a8}
\ee
By (\ref{a7}) with $I$ replaced by $K$, $J_{C}(k_1) = I_A(j)$,
and (\ref{a5}) follows. It remains to prove
(\ref{a8}), and we shall use monotonicity again for this.

It suffices for (\ref{a8}) to prove that
\be
\la(C\mid U_j=1) \ge \la(B\mid U_j=1),
\label{a9}
\ee
together with the reversed inequality given $U_j=0$.
The conditioning on the left-hand side of (\ref{a9})
refers to the first coordinate encountered by the
algorithm (\ref{a1})--(\ref{a2'}) when applied to the re-ordered
set $K$.
Let 
\be
\ol U=
(U_1,U_2,\dots,U_{j-1},1,U_{j+1},\dots, U_N).
\label{old1}
\ee
The $0/1$-vector $f(\ol U)=(f_i(\ol U): 1\le i\le N)$ 
is constructed sequentially (as above) by considering the indices
$1,2,\dots,N$ in turn. At stage $k$, we declare
$f_k(\ol U)$ to equal 1 if $U_k$ exceeds a certain
function $a_k$ of the variables $f_i(\ol U)$, $1\le i < k$.
By the monotonicity of $\mu$, this
function is non-increasing in these 
variables.
The index $j$ plays a special role in that: (i) $f_j(\ol U) = 1$,
and (ii) given this fact, it is more likely than before
that the variables $f_k(\ol U)$, $j<k\le N$, will
take the value 1. The values $f_k(\ol U)$, $1\le k<j$ are unaffected
by the value of $U_j$.

Consider now the $0/1$-vector $g(\ol U) =
(g_{k_r}(\ol U): 1\le r\le N)$, constructed in the same manner
as above but with the new ordering $K$ of the index set $I$.
First we examine index $k_1$ ($=j$), and we 
automatically declare $g_{k_1}(\ol U)
=1$ (since $U_j=1$). We then construct $g_{k_r}(\ol U)$, 
$2\le r\le N$,
in sequence. Since the $a_k$ are non-decreasing
in the variables constructed so far, we have that
\be
g_{k_r}(\ol U) \ge f_{k_r}(\ol U),\qq
r=2,3,\dots,N.
\label{new23}
\ee
Therefore, $g(\ol U) \ge f(\ol U)$,
implying as required that
\be
\la(C\mid U_j=1) =\la (g(\ol U) \in A) \ge \la(f(\ol U)\in A)
=\la(B\mid U_j=1).
\label{a11}
\ee
Inequality (\ref{a9}) follows. The same argument implies the 
reversed inequality
obtained from (\ref{a9}) by reversing the conditioning to $U_j=0$.
This implies (\ref{a8}).

A formal proof of (\ref{new23}) follows.
Suppose that $r$ is such that $g_{k_s}(\ol U) \ge f_{k_s}(\ol U)$
for $2\le s < r$.  By (\ref{a2'}), for $r\le j$,
\begin{align*}
f_{k_r}(\ol U) &= 1\q\text{if} \q U_{k_r} > \mu(X_{k_r}=0\mid
   X_{k_s}=f_{k_s}(\ol U)\text{ for } 2\le s < r),\\
g_{k_r}(\ol U) &= 1\q\text{if} \q U_{k_r} > \mu(X_{k_r}=0\mid
   X_{k_s}=g_{k_s}(\ol U)\text{ for } 1\le s < r).
\end{align*}
Now $g_{k_1}(\ol U) =1$ and, by the induction hypothesis and
monotonicity,
\begin{multline*}
\mu(X_{k_r}=0\mid
 X_{k_s}=f_{k_s}(\ol U)\text{ for } 2\le s < r)\\
 \ge \mu(X_{k_r}=0\mid
      X_{k_s}=g_{k_s}(\ol U)\text{ for } 1\le s < r),
\end{multline*}
whence $g_{k_r}(\ol U) \ge f_{k_r}(\ol U)$ as required.

Consider finally the case $j < r\le N$. Then
\begin{align*}
f_{k_r}(\ol U) &= 1\q\text{if} \q U_{k_r} > \mu(X_{k_r}=0\mid
   X_{k_s}=f_{k_s}(\ol U)\text{ for } 1\le s < r),\\
   g_{k_r}(\ol U) &= 1\q\text{if} \q U_{k_r} > \mu(X_{k_r}=0\mid
      X_{k_s}=g_{k_s}(\ol U)\text{ for } 1\le s < r),
\end{align*}
and the conclusion follows as before.
\end{proof}

\section{Sharp-threshold theorem}\label{sectst}
We consider in this section a family of probability
measures indexed by a parameter
$p\in(0,1)$, and we prove a sharp-threshold theorem 
subject to a hypothesis of monotonicity. The motivating
example is the random-cluster model, to which we return in the next section.

Let $1\le N<\oo$, $I=\{1,2,\dots,N\}$,
and let $\Om=\{0,1\}^N$ and $\sF$ be given
as before.
Let $\mu$ be a positive probability measure on $(\Om,\sF)$.
For $p\in (0,1)$, we define the probability measure
$\mu_p$ by
\be
\mu_p(\om) = \frac 1{Z_p} \mu(\om)
\left\{\prod_{i\in I} p^{\om(i)}(1-p)^{1-\om(i)}\right\},
\qq \om\in\Om,
\label{mupdef}
\ee
where $Z_p$ is the normalizing constant
\be
Z_p = \sum_{\om\in\Om}\mu(\om)
\left\{\prod_{i\in I} p^{\om(i)}(1-p)^{1-\om(i)}\right\}.
\ee
It is immediate that $\mu_p$ is positive and that $\mu=\mu_{\frac12}$.
It is 
easy to check that $\mu_p$ satisfies the FKG lattice condition
(\ref{FKG}) if and only if $\mu$ satisfies this condition,
and it follows that $\mu$ is monotonic if and only if, for all
$p\in (0,1)$, $\mu_p$ is monotonic.
In order to prove a sharp-threshold theorem for the family
$\mu_p$, we present first a Russo-type formula.

\begin{thm}[\cite{BGK}]\label{russo}
For any event $A\in\sF$,
\be
\frac{d}{dp}\mu_p(A) = \frac 1{p(1-p)} 
\sum_{i\in I} \cov_p(X_i,1_A),
\label{russodiff}
\ee
where $\cov_p$ denotes covariance with 
respect to the measure $\mu_p$.
\end{thm}

\begin{proof} 
This may be obtained exactly as in \cite{BGK}, Proposition 4,
see also Section 2.4 of \cite{G-RC}. The details are omitted.
\end{proof}

Let $\sA$ be a subgroup of the permutation group
$\Pi_N$. A probability measure $\phi$ on $(\Om,\sF)$ is
called {\it $\sA$-invariant\/} if $\phi(\om)=\phi(\a\om)$ for
all $\a\in\sA$. An event $A\in\sF$ is
called {\it $\sA$-invariant\/} if $A=\a A$ for
all $\a\in\sA$. It is easily seen that, for any subgroup $\sA$,
$\mu$ is $\sA$-invariant if and only if each $\mu_p$ is
$\sA$-invariant.

\begin{thm}[Sharp threshold]\label{genconc}
There exists a constant $c$ satisfying $c\in(0,\oo)$
such that the following holds. Let $N\ge 1$ and let $A\in\sF$
be an increasing 
event. 
Let $\mu$ be a positive probability
measure on $(\Om,\sF)$ which is monotonic. If there
exists a subgroup $\sA$ of $\Pi_N$ acting transitively
on $I$ such that $\mu$ and $A$ are $\sA$-invariant, then
\be
\frac{d}{dp}\mu_p(A) \ge
  \frac{c\xi_p}{p(1-p)}
\min\{\mu_p(A), 1-\mu_p(A)\} \log N,\qq p\in(0,1),
\label{shpt}
\ee
where $\xi_p =\min\{\mu_p(X_i)(1-\mu_p(X_i)): i\in I\}$.
\end{thm}

We precede the proof with a lemma. Let 
$$
I_{p,A}(i) = \mu_p(A \mid X_i=1) - \mu_p(A \mid X_i=0).
$$

\begin{lem}\label{lemma}
Let $A\in\sF$. Suppose there
exists a subgroup $\sA$ of\/ $\Pi_N$ acting transitively
on $I$ such that $\mu$ and $A$ are $\sA$-invariant.
Then $I_{p,A}(i)=I_{p,A}(j)$ for all $i,j\in I$ and all $p\in (0,1)$.
\end{lem}

\begin{proof}[Proof of Lemma \ref{lemma}]
Since $\mu$ is $\sA$-invariant, so is
$\mu_p$ for every $p$. 
Let $i,j\in I$, and find $\a\in\sA$ such that $\a_i=j$.
Under the given conditions,
\begin{align*}
\mu_p(A,\, X_j=1) &= \sum_{\om\in A}\mu_p(\om)X_j(\om)
=\sum_{\om\in A} \mu_p(\a\om)X_{i}(\a\om)\\
&=\sum_{\om'\in A} \mu_p(\om')X_i(\om')
= \mu_p(A,\, X_i=1).
\end{align*}
Applying this with $A=\Om$, we find that
$\mu_p(X_j=1)=\mu_p(X_i=1)$. By dividing, we deduce
that $\mu_p(A\mid X_j=1)=\mu_p(A\mid X_i=1)$.
A similar equality holds with 1 replaced by 0, and the 
claim follows.
\end{proof}

\begin{proof}[Proof of Theorem \ref{genconc}]
By Lemma \ref{lemma},
every index has the same influence. 
Since $A$ is increasing,
\begin{align*}
\cov_p(X_i,1_A) &= \mu_p(X_i1_A) - \mu_p(X_i)\mu_p(A)\\
&=  \mu_p(X_i)(1-\mu_p(X_i)) I_{p,A}(i)\\
&\ge \xi_p I_{p,A}(i).
\end{align*}
Summing over the index set $I$ as in 
(\ref{russodiff}), we deduce (\ref{shpt}) by Theorem 
\ref{infthm} applied to the monotonic measure $\mu_p$.
\end{proof}

\section{Probability measures on the Euclidean cube}
We have so far considered probability measures on
the discrete cube $\{0,1\}^N$ only.  The method of proof
of the influence theorem, Theorem \ref{infthm}, may be applied also
to probability measures on the Euclidean cube $[0,1]^N$ that are
absolutely continuous with respect to Lebesgue measure. Any
such measure $\mu$ has a density function $\f$, which is to say that
$$
\mu(A) = \int_A\f(\bx)\,\lambda(d\bx),
$$
for (Lebesgue)
measurable subsets $A$ of $[0,1]^N$, with $\lambda$ denoting
Lebesgue measure. Since the density function $\f$ is non-unique,
we shall phrase the results of this section in terms of
$\f$ rather than the associated measure $\mu$. Some may regard
this as not entirely satisfactory, arguing that results for
{\it measures\/} should be based on hypotheses for these measures,
rather than for particular versions of their density functions.
One may rewrite the conclusions of this section thus, but 
at the expense of greater measure-theoretic detail which obscures the
basic argument.

Let $N\ge 1$, and write $\Om=[0,1]^N$.
Let $\f:\Om\to[0,\oo)$ be (Lebesgue) measurable. We call $\f$
a {\it density function\/} if
$$
\int_{\Om}\f(\bx)\,\la(d\bx) = 1,
$$
and in this case we denote by $\mu_\f$ the corresponding probability
measure,
$$
\mu_\f(A) = \int_A \f(\bx)\,\la(d\bx).
$$
We call
$\f$ {\it positive\/} if it
is a strictly positive function on $\Om$, and we say
it satisfies the {\it (continuous) FKG lattice
condition\/} if
\be
\f(\bx\vee \by)\f(\bx\wedge \by) \ge \f(\bx)\f(\by)\qq
\text{for all } \bx,\by\in \Om,
\label{new24}
\ee
where the operations $\vee$, $\wedge$ are defined as the 
coordinate-wise maximum
and minimum, respectively.

Let $\f$ be a density function. 
We call $\mu_\f$ {\it positively associated\/}
if
$$
\mu_\f(A\cap B) \ge \mu_\f(A)\mu_\f(B),
$$
for all increasing subsets of $\Om$.
[It is presumably well known that increasing subsets of $\Om$ are
Lebesgue-measurable but need not be Borel-measurable; see the notes
at the end of this section.] 

Let $I=\{1,2,\dots,N\}$.
For $J\subseteq I$, let
$\Om_J=[0,1]^J$ and
\be
\Om_J^\xi=\{\bx\in\Om:x_j=\xi_j \text{ for } j \in I\sm J\},
\qq \xi\in\Om.
\label{new4.2}
\ee
The Lebesgue $\s$-algebra of $\Om_J$ is denoted by $\sF_J$.
Let $\f$ be a positive density function. We define the
conditional probability measure $\mu_{\f,J}^\xi$ on $(\Om_J,\sF_J)$  by
\be
\mu_{\f,J}^\xi(E) = \int_{E} \f_{J}^\xi(\bx)
\,\la(d(x_j:j\in J)),
\qq E\in \sF_J,
\label{new4.3}
\ee
where $\f_J^\xi$ is the conditional density function
$$
\f_{J}^\xi(\bx)
 = \frac1{Z_J^\xi} \f(\bx)1_{\Om_J^\xi}(\bx),\q
Z_J^\xi =  \int_{\Om_J^\xi}\f(\bx)\,\la(d(x_j:j\in J)).
$$
We sometimes write $\mu_\f\bigl(E\mid(\xi_j: j\in I\sm J)\bigr)$
for $\mu_{\f,J}^\xi(E)$, and we recall the standard fact that
$\mu_\f\bigl(\cdot\mid (\xi_j:j\in I\sm J)\bigr)$
is a version of the conditional expectation
given the $\s$-field $\sF_{I\sm J}$.

We say that $\f$
is {\it strongly positively-associated\/}
if:
for all $J\subseteq I$ and all $\xi\in\Om$, the measure
$\mu_{\f,J}^\xi$ is positively associated.
We call $\f$
{\it monotonic\/} if: for all $J\subseteq I$, all
increasing subsets $A$ of $\Om_J$, and
all $\xi,\zeta\in\Om$,
\be
\mu_{\f,J}^\xi(A) \le
\mu_{\f,J}^\zeta(A)\qq\text{whenever } \xi\le\zeta.
\label{new4.3a}
\ee
That is, $\f$ is monotonic if, for all $J\subseteq I$, 
\be
\mu_{\f,J}^\xi\lest\mu_{\f,J}^\zeta \qq\text{whenever } \xi\le\zeta.
\label{new4.4} 
\ee

Here is a basic result concerning stochastic ordering.

\begin{thm} [\cite{BattyBollmann, Preston}]\label{BBP}
Let $N\ge 1$, and let $f$ and $g$ be density functions on 
$\Om=[0,1]^N$.
If
\[
g(\bx\vee \by)f(\bx\wedge \by) \ge g(\bx)f(\by)\qq
 \text{for all } \bx,\by\in [0,1]^N,
\]
then $\mu_f \lest \mu_g$.
\end{thm}

If $\f$ satisfies the FKG lattice condition and $A$ is an increasing event,
then
$$
1_A(\bx\vee\by)\f(\bx\vee \by)\f(\bx\wedge \by) \ge 
1_A(\bx)\f(\bx)\f(\by),
$$
whence, by Theorem \ref{BBP},
$$
\mu_\f(A)\mu_\f(B) \le \mu_\f(A \cap B)
$$
for all increasing $A$, $B$. Therefore, $\mu_\f$ is positively
associated. 

Henceforth we restrict ourselves to {\it positive\/}
density functions. Arguments similar to the above
are valid with $\f$ (assumed positive) replaced
by the conditional density function $\f_J^\xi$, and one arrives 
thus at the following.

\begin{thm}\label{ctspa}
Let $N\ge 1$, and let $\f$ be a positive density function
on $\Om=[0,1]^N$ satisfying the FKG lattice condition (\ref{new24}).
Then $\f$ is strongly positively-associated and monotonic.
\end{thm}

We turn now to a `continuous' version of Theorem \ref{infthm}.
Let $N\ge 1$, and let
$\f$ be a monotonic positive density function
on $\Om=[0,1]^N$. 
Let $U=(U_1,U_2,\dots,U_N)$ be the identity function
on $[0,1]^N$.
For an increasing subset $A$ of $\Om$,
we define the {\it\coin{}s\/} by
\be
I_A(i) = \mu_\f(A\mid U_i=1) - \mu_\f(A\mid U_i=0),\qq i\in I.
\label{infcontinuous}
\ee

\begin{thm}[Influence]\label{infthm2}
There exists a constant $c$ satisfying $c\in(0,\oo)$
such that the following holds. Let $N\ge 1$ and let $A$
be an increasing
subset of $\Om=[0,1]^N$.
Let $\f$ be a positive density function on $[0,1]^N$
that is monotonic.
There exists $i\in I$ such that
\be
I_A(i) \ge c\min\{\mu(A),1-\mu(A)\} \frac{\log N}N.
\ee
\end{thm}

\begin{proof}
The proof is very similar to that of Theorem \ref{infthm}.
We propose first to construct an increasing 
event $B$ such that $\la(B)=\mu(A)$,
by way of a function $f:[0,1]^N\to[0,1]^N$.
Let  $\bx=(x_i:1\le i\le N)\in[0,1]^N$, and write
$f(\bx)=(f_1(\bx),f_2(\bx),\dots,f_N(\bx))$.
The first coordinate $f_1(\bx)$ depends on $x_1$ only
and is defined by:
\[
\mu_\f(U_1>f_1(\bx))=1-x_1.
\]
Since the density function $\f$ is strictly
positive, $f_1(\bx)$ is a continuous and strictly increasing function of
$x_1$. It is an elementary exercise to check
that the law of $f_1(U)$ under $\la$ is the same as
that of $U_1$ under $\mu_\f$.

Having defined $f_1(\bx)$, we define $f_2(\bx)$ in terms of $x_1$, $x_2$
only by:
\[
\mu_\f\bigl(U_2>f_2(\bx)\bigmid U_1=f_1(\bx)\bigr)
=1-x_2.
\]
The left-hand side is defined according to (\ref{new4.3}). 
It is a standard fact that
$\mu_\f(\cdot\mid U_1=f_1)$ is a version of the conditional
expectation $\mu_\f(\cdot\mid \s(U_1))$, where $\s(U_1)$ denotes the $\s$-field
generated by $U_1$, and it is an exercise to check that the pair
$(f_1(U),f_2(U))$ has the same law under $\la$ as
does the pair $(U_1, U_2)$ under $\mu_\f$. For each given $x_1\in(0,1)$,
$f(\bx)$ is a continuous and strictly increasing function of $x_2$.
[We use the assumptions that $\f$ is positive and monotonic,
respectively, here.]

We continue inductively.
Suppose we know $f_i(\bx)$ for $1\le i < k$. Then $f_k(\bx)$ 
depends on $x_1,x_2,\dots,x_k$ and is
given by:
\[
\mu_\f\bigl(U_k>f_k(\bx)\bigmid U_i=f_i(\bx)\text{ for } 1\le i<k\bigr)
=1-x_k.
\]
As above, $f$ is strictly increasing (using the assumption of
monotonicity), and the law of 
$f(U)$ under $\la$ is the same
as the law of $U$ under $\mu_\f$. We set $B=f^{-1}(A)$.

Let
$$
J_B(i) = \la(B\mid U_i=1) - \la(B\mid U_i=0),\qq i\in I.
$$
Since $f_1$ is continuous and strictly increasing, 
$$
\mu_\f(A\mid U_1=b) = \la(B\mid f_1(U_1)=b)=\la(B \mid U_1=b),
\qq b=0,1,
$$
implying that $I_A(1)=J_B(1)$.
It remains to show that $I_A(j) \ge J_B(j)$ for $j\in I$.
Let $j\in I$, $j\ne 1$.
We re-order the coordinate set as 
$K=\{j,1,2,\dots,j-1,j+1,\dots,N\}$,
and we construct a continuous
increasing function $g$ as above but subject to the new ordering.
Rather than re-work the details from the
proof of Theorem \ref{infthm}, we prove only part of that necessary.
We sketch a proof that $\mu_\rho(A\mid U_j=1)
\ge \la(B\mid U_j=1)$, a similar argument being valid with 1
replaced by 0 and the inequality reversed.
The main step is to show that $f\le g$ under the assumption that $U_j=1$.
Suppose that $1\le r<j$, and assume it has already been proved
that $f_i(\bx)\le g_i(\bx)$ for $\bx\in\Om$
and $1\le i<r$.
Let $\bx\in\Om$. We claim that
\begin{multline}
\mu_\f(U_r>\xi\mid U_i=f_i(\bx)\text{ for } 1\le i<r)\\
\le
\mu_\f(U_r>\xi\mid U_j=1,\  U_i=g_i(\bx)\text{ for } 1\le i<r),
\qq \xi\in [0,1].
\label{new5.1}
\end{multline}
By monotonicity,
\begin{multline}
\mu_{\rho,J}(\cdot\mid U_j=u,\ U_i=f_i(\bx)\text{ for } 1\le i<r)
\\
\lest  \mu_{\rho,J}(\cdot\mid 
    U_j=1,\ U_i=g_i(\bx)\text{ for } 1\le i<r),\qq u\in[0,1].
\label{new5.2}
\end{multline}
The left-hand side of (\ref{new5.2}) is a version of the conditional
expectation of the conditional measure $\mu_{\rho,J}(\cdot\mid U_i=f_i(\bx)\text{ for } 1\le i<r)$ given $\s(U_j)$. By averaging
over the value of $u$ in (\ref{new5.2}), we obtain (\ref{new5.1}).
The other steps are proved similarly.
\end{proof}
 
Unlike the discrete setting of Section 3,
Theorem \ref{infthm2} does not imply a sharp-threshold theorem.
Any density function $\f$ on
$[0,1]^N$ may be used to generate a parametric family
$(\f_p: 0< p< 1)$ of densities given by
$$
\f_p(\bx) = \frac 1{Z_{\f,p}}\f(\bx)\prod_{i=1}^N p^{x_i}(1-p)^{1-x_i},
\qq \bx=(x_1,x_2,\dots,x_N)\in[0,1]^N,
$$
and we write $\mu_p=\mu_{\f_p}$.
Let $A$ be an increasing
subset of $[0,1]^N$. The proof of Theorem \ref{russo} may be adapted to
this setting to obtain that
$$
\frac d{dp} \mu_{p}(A)=\frac 1{p(1-p)}\sum_{i=1}^N \cov_p(U_i,1_A),
$$
where $U=(U_1,U_2,\dots,U_N)$ is the identity function on $[0,1]^N$,
and
$\cov_p$ denotes covariance with respect to $\mu_p$.

Let $\f$ be the constant function, so that $\mu_\f$ is Lebesgue
measure. As above,
let $p\in (0,1)$ and let
$Y_1,Y_2,\dots,Y_N$ be independent random variables
taking values in $[0,1]$ with common density function
$$ 
\f_p(x) =  
\begin{cases}\dfrac{\log[p/(1-p)]}{2p-1} p^x(1-p)^{1-x}
  \q&\text{if } p\ne \frac12,\ x\in (0,1),\\
  1 &\text{if } p=\frac12,\ x\in (0,1).
  \end{cases}
$$
It is easily checked that the joint density function
$$
\f_p(\bx) = \prod_{i=1}^N \f_p(x_i), \qq \bx=(x_1,x_2,\dots,x_N)\in[0,1]^N,
$$
satisfies the FKG lattice condition, and is therefore monotonic.

We now choose $A$ by $A= (N^{-1},1]^N$. It is an easy calculation
that
$$
\mu_p(A) =\begin{cases} \left(1-\dfrac{\pi^{1/N}-1}{\pi-1}\right)^N
 \q&\text{if } p\ne\frac12,\\
 \left(1-\dfrac1N\right)^N&\text{if } p=\frac12,
\end{cases}
$$
where $\pi=p/(1-p)$.
Therefore, as $N\to\oo$,
$$
\mu_p(A) \to \begin{cases} \pi^{-1/(\pi-1)}\q&\text{if } p\ne \frac12,\\
e^{-1} &\text{if } p=\frac12.
\end{cases}
$$
In addition, 
$$
\cov_{\frac12}(U_i,1_A) = \frac 1N\left(1-\frac1N\right)^{N-1}
\sim \frac{e^{-1}}N.
$$
The influence theorem, Theorem \ref{infthm2}, may be applied
to the event $A$, but there is no sharp threshold for $\mu_p(A)$.
This situation diverges from that of the discrete setting
at the point where a lower bound for the \coin\ $I_A(i)$ is
used to calculate a lower bound for the covariance $\cov_p(U_i,1_A)$.
    
We return briefly to the measurability of an increasing subset
of $[0,1]^N$.

\begin{thm}\label{Lmeas}
Let $N\ge 2$. Every
increasing subset of $[0,1]^N$ is Lebesgue-measurable.
\end{thm}

Increasing subsets need not be Borel-measurable, as the following
example indicates. Let $M$ be a non-Borel-measurable
subset of $[0,1]$.  Consider the increasing subset $A$ of $[0,1]^2$
given by
$$
A= \{(x,y)\in[0,1]^2: x+y>1\} \cup\{(x,1-x): x\in M\}.
$$
The function $h: x\mapsto (x,1-x)$ is a continuous, and hence
Borel-measurable, function
from $\RR$ to $\RR^2$. If $A$ were Borel-measurable, then so would be
$$
A'=A\cap \{(x,1-x): x\in \RR\} = \{(x,1-x): x\in M\}.
$$ 
This would imply that $h^{-1}(A') = M$ is Borel-measurable, a contradiction.

\begin{proof}[Proof of Theorem \ref{Lmeas}]
The statement is trivially true when $N=1$, and we prove the general
case by induction on $N$. Suppose $n$ is such
that the result holds for $N=n$.
Let $A$ be an increasing subset of $[0,1]^{n+1}$, and let 
$g:[0,1]^n\to[0,1]\cup\{\oo\}$ be defined by
\[
g(\bx)=\inf \{y: (\bx,y) \in A\},\qq \bx\in[0,1]^n.
\]
The function $g$ is decreasing on $[0,1]^n$, and hence, for all $c\in\RR$,
the subset $H_c=\{\bx:g(\bx)<c\}$ is increasing. By the induction hypothesis,
each $H_c$ is Lebesgue-measurable in $[0,1]^n$, and therefore
$g$ is a measurable function. 
Its graph
$G=\{(\bx,g(\bx)): \bx\in[0,1]^{n}\}$ is (by an approximation
by simple functions, or otherwise) a Lebesgue-measurable set
and is also (by Fubini's Theorem) a null
subset of $[0,1]^{n+1}$. Furthermore, the set
$$
\ol A = \{(\bx,y)\in [0,1]^{n+1}: y > g(\bx)\}
$$ 
is
Lebesgue-measurable.
Now $A$ differs from $\ol A$ only on a subset of the null set $G$,
and the claim follows.
\end{proof}

\section{The random-cluster model}
The sharp-threshold theorem of Section \ref{sectst} may 
be applied as follows to
the \rc\ measure.
Let $G=(V,E)$ be a finite graph, assumed for simplicity
to have neither loops nor multiple edges. 
We take as configuration space the set 
$\Om=\{0,1\}^E$, and write $\sF$ for the set of its subsets.
For $\om\in\Om$, we call
an edge $e$ {\it open\/} (in $\om$) if $\om(e)=1$, 
and {\it closed\/} otherwise. Let $\eta(\om)=
\{e\in E: \om(e)=1\}$ be the set of
open edges, and consider the open graph $G_\om=(V,\eta(\om))$.
The connected components of $G_\om$ are termed
{\it open clusters\/}, and $k(\om)$ denotes
the number of such clusters (including 
any isolated vertices).

Let $q\in(0,\oo)$, and let $\mu$ be the probability measure
on $(\Om,\sF)$ given by
\be
\mu(\om) = \frac 1 {Z(q)} 
q^{k(\om)}, \qq \om\in\Om,
\label{defmu}
\ee
where $Z(q)$ is the appropriate normalizing constant.
It is clear that $\mu$ is positive, and it is easily checked
that $\mu$ satisfies the FKG lattice condition 
if $q\ge 1$. See
\cite{F72b, G-RC}.
(The FKG lattice condition does not hold
when $q<1$ and $G$ contains a circuit.)
{\it We assume henceforth
that $q\ge 1$.} 
By Theorem \ref{new2.7},
$\mu$ is monotonic.

The \rc\ measure $\fpq$ on the graph $G$
with parameters $p\in(0,1)$ and $q\in[1,\oo)$
is given as in (\ref{mupdef}) by
\be
\fpq(\om) = \frac 1{Z(p,q)} \left\{
\prod_{e\in E} p^{\om(e)}(1-p)^{1-\om(e)}\right\}
q^{k(\om)},
\qq \om\in\Om.
\label{mupdef2}
\ee
It is well known (see \cite{F72b, G-RC}) that
\be
\frac p{p+q(1-p)} \le \fpq(X_e=1) \le p,\qq e\in E.
\label{finen}
\ee
We call $G$ $\sA$-{\it transitive\/} if its automorphism
group possesses a subgroup $\sA$
acting transitively on $E$. 
We may apply Theorem \ref{genconc} to obtain the
following.
There exists an absolute
constant $c>0$ such that, for all $\sA$-transitive graphs $G$,
all $p$, $q$,  and 
any increasing $\sA$-invariant event $A\in\sF$,
\be
\frac{d}{dp}\fpq(A) \ge
 c \min\left\{\frac{q}{\{p+q(1-p)\}^2}, 1\right\} 
  \min\{\fpq(A), 1-\fpq(A)\} \log N,
\nonumber
\ee
whence
\be
\frac{d}{dp}\fpq(A) \ge
  \frac cq \min\{\fpq(A), 1-\fpq(A)\} \log N.
\label{fpqsteep}
\ee
The differential inequality (\ref{fpqsteep}) takes the usual
simpler form when $q=1$, and it 
may be integrated
exactly for general $q\ge 1$.
Here is an illustration
of (\ref{fpqsteep}) when integrated. Let $p_1\in(0,1)$
be chosen such that $\phi_{p_1,q}(A)\ge\frac12$, and let $p_1<p_2<1$.
We note that $\fpq(A) \ge \frac12$
for $p\in (p_1,p_2)$.
We integrate (\ref{fpqsteep})
over this interval to obtain that
\be
\phi_{p_2,q}(A) \ge 1-\tfrac12 N^{-c(p_2-p_1)/q}.
\label{upper}
\ee

Bollob\'as and Riordan have shown in \cite{BR2, BR1} how to apply
the sharp-threshold theorem for product measure to percolation in
two dimensions, thereby obtaining a further proof of the
famous theorem of Harris and Kesten that the critical probability of
bond percolation equals $\frac12$. Their key step
is the proof that there exists a sharp threshold for the 
event that a large square is traversed by an open path.
One obtains similarly the following for
the \rc\ model on the
square lattice $\LL^2$.

Let $\ZZ=\{\dots,-1,0,-1,\dots\}$ be the integers,
and $\ZZ^2$ the set of all $2$-vectors $x=(x_1,x_2)$
of integers. We turn $\ZZ^2$ into a graph by placing an edge
between any two vertices $x$, $y$ with $|x-y|=1$,
where
$$
|z| = |z_1| + |z_2|, \qq z\in \ZZ^2.
$$
We write $\EE^2$ for the set of such edges, and $\LL^2=(\ZZ^2,\EE^2)$
for the ensuing graph. We shall work on a finite torus of $\LL^2$.
Let $n\ge 1$. Consider the square $S_n=[0,n]^2$ (this is a convenient
abbreviation for $\{0,1,2,\dots,n\}^2$) viewed as a subgraph
of $\LL^2$. We identify certain pairs of
 vertices on the boundary of $S_n$ 
in order to make it symmetric. More specifically, we identify
any pair of the form $(0,m)$, $(n,m)$ and of
the form $(m,0)$,
$(m,n)$, for $0\le m\le n$, and we merge any parallel
edges that ensue. Let $T_n=(V_n,E_n)$ denote the 
resulting toroidal graph. 
Let $\sA_n$ be the automorphism group of the graph $T_n$, and note
that $\sA_n$ acts transitively on $E_n$. The configuration space of
the \rc\ model on $T_n$ is denoted $\Om(n)=\{0,1\}^{E_n}$.

Let $p\in(0,1)$ and
$q\in[1,\oo)$. Write $\phi_{n,p}$
for the \rc\ measure on $T_n$ with parameters $p$ and $q$,
and note that $\phi_{n,p}$ is $\sA_n$-invariant.
Let 
$$
\psd=\psd(q) = \frac{\sqrt q}{1+\sqrt q},
$$
the self-dual point of the \rc\ model on $\LL^2$, see
\cite{G03, G-RC}. 
We note that the (Whitney)
dual of $T_n$ is isomorphic to $T_n$, and the \rc\ measure 
on $T_n$ is self-dual
when $p=\psd$.

Let $\om\in\Om(n)$.
Any translate in $T_n$ of a rectangle of the form
$[0,r]\times[0,s]$ is said to be
of size $r\times s$. When $r\ne s$, such a translate
is said to be traversed {\it long-ways\/} 
(\resp, traversed {\it short-ways\/})  if the two shorter
sides  (\resp, longer sides) of the rectangle are
joined within the rectangle by an open path of $\om$.

Let $k\ge 2$, $n\ge 1$.
Let $R_n=[0,n+1]\times[0,n]$, viewed as a subgraph of $T_{kn}$,
and let $\LW_n$ be the event that $R_n$ is traversed long-ways.
By a standard duality argument,
\be
\phi_{kn,\psd}(\LW_n)=\tfrac12,\qq k\ge 2,\ n\ge 1.
\label{selfdual}
\ee
Let $A_n$ be the event that there exists in $T_{kn}$ some
translate of the square $S_n=[0,n]\times[0,n]$ that
possesses either an open top--bottom crossing or
an open left--right crossing. The event
$A_n$ is $\sA_n$-invariant, and
\be
\phi_{kn,\psd}(A_n)\ge \phi_{kn,\psd}(\LW_n) = \tfrac12.
\label{greater}
\ee

We apply (\ref{upper}) to the event
$A_n$, with $p_1=\psd$ and with $N=2(kn)^2$ being the number of edges
in $T_{kn}$. This yields 
that
\begin{align}
\phi_{kn,p}(A_n) &\ge 1 - \tfrac12 [2(kn)^2]^{-c(p-\psd)/q}
\nonumber\\
&\ge 1 - (kn)^{-2c(p-\psd)/q},
\qq \psd<p <1.
\label{upper4}
\end{align}

The event $A_n$ is defined on the whole of the torus. We next
use an argument taken from \cite{BR2, BR1} to obtain a more
locally defined event.  We shall for simplicity of notation
treat certain
real-valued quantities as if they were integers. Let $1<\a <k$,
and let $H_{n,\a} = [0,\a n]\times[0,n/\a]$ and
$V_{n,\a} = [0,n/\a]\times[0,\a n]$. 
Let $h_{n,\a}$, $v_{n,\a}$ be the sets of vertices
in $T_{kn}$ given by
\begin{align*}
h_{n,\a} &= \{ (l_1n(\a-1), l_2n(1-\a^{-1}))\in V_{kn}: l_1,l_2\in \ZZ\},\\
v_{n,\a} &= \{(l_1n(1-\a^{-1}), l_2n(\a -1))\in V_{kn}: l_1,l_2\in \ZZ\}.
\end{align*}
Consider the set $\sH=H_{n,\a} + h_{n,\a}$ of translates of $H_{n,\a}$ by
vectors in $h_{n,\a}$, and also the set $\sV=V_{n,\a}+v_{n,\a}$. If $A_n$
occurs, then some rectangle in $\sH\cup\sV$ is traversed
short-ways. By positive association and symmetry,
\begin{align}
\phi_{kn,p}(\comp{A_n}) &\ge
\phi_{kn,p}(\text{no member of $\sH\cup\sV$ is traversed short-ways})
\nonumber\\
&\ge \{1-\phi_{kn,p}(\SW_{n,\a})\}^M,
\end{align}
\label{upper5}
where $\SW_{n,\a}$ is the event that $H_n$ is traversed short-ways,
and 
\be
M=|h_{n,\a}| + |v_{n,\a}|.
\label{defM}
\ee
After taking into account the rounding effects above,
we find that
\be
M \le 2\left(1+\frac{k}{\a-1-n^{-1}}\right)
\left(1+\frac{k}{1-\a^{-1}-n^{-1}}\right),
\label{Mlower}
\ee
so that $M$ is approximately $2k^2\a/(\a-1)^2$
when $k$ and $n$ are large.

Combining (\ref{upper4})--(\ref{defM}), we arrive at the
following theorem, where $\SW_{n,\a}$ is the event
that the rectangle $\bigl[0,\lfloor n\a\rfloor\bigr]\times
\bigl[0,\lfloor n/\a\rfloor\bigr]$ is crossed short-ways.

\begin{thm}\label{thmexcat} Let $k\ge 2$, $n\ge 1$, and
$\psd< p < 1$. We have that
\be
\phi_{kn,p}(\SW_{n,\a})
\ge 1- e^{- g(p-\psd)}
\label{upper6}
\ee
where
$$
g=g(k,n,\a,q) = \frac {2c}{Mq} \log(kn).
$$
\end{thm}

In particular, for $p>\psd$, one may make $\phi_{kn,p}(\SW_{n,\a})$
large by holding $k$ fixed and sending $n\to\oo$.
It does not seem to be easy to deduce an estimate for $\fpq(\SW_{n,\a})$
for a \rc\ measure $\fpq$ on the infinite lattice $\LL^2$.
Neither do we know how to use the existence of crossings short-ways
to build crossings long-ways. This is in contrast
to the case of product measure, see \cite{BR1, CC, G99, Ru78, Ru81, SeW}.

\section{The critical point}
There is a famous conjecture that the critical point $\pc(q)$ of
the \rc\ model on $\LL^2$ equals $\psd (q)$. We do not
spell out the details necessary to state this conjecture
properly, referring the reader instead to \cite{G03, G-RC}.
The conjecture is known to be valid for $q=1$ (percolation),
$q=2$ (a case corresponding to the Ising model), and for sufficiently
large $q$ (namely $q \ge 21.61$).  The conjecture would follow
if one could prove a strengthening of Theorem \ref{thmexcat} in which
short-ways is replaced by long-ways, and with the toroidal
measure replaced by the wired measure on the full lattice.
We finish by explaining this. 

The so-called
`wired \rc\ measure' on $\LL^2$ is denoted by
$\fpq^1$, and the reader is referred to the references
above for a definition of $\fpq^1$. 

\begin{thm}\label{alternativefinaltheorem}
Let $q\ge 1$.
Let $p_k$ be the $\fpq^1$-probability that a 
$2^k \times  2^{k+1}$ 
rectangle is crossed long-ways.
Suppose that 
\be
\prod_{k=1}^{\oo}  p_k > 0, \qq p > \psd(q).
\label{finalcond}
\ee
Then the critical point of the \rc\ model on $\LL^2$ equals $\psd(q)$.
\end{thm}

By duality,
$1-p_k = \phi_{p',q}^0(\SW(k))$,
where $\SW(k)$ is the event that the rectangle
$[0,2^{k+1}-1]\times[0,2^k +1]$
is traversed short-ways, and $p'$ is the dual value of $p$,
$$
\frac {p'}{1-p'} = \frac{q(1-p)}p.
$$
Therefore,
\begin{align*}
\sum_{k=1}^\oo(1-p_k) &\le
\sum_{k=1}^\oo 2^{k+1}\phi_{p',q}^0(\rad(C)\ge 2^k+1)\\
&\le 4\sum_{n=1}^\oo \phi_{p',q}^0(\rad(C)\ge n)\\
&= 4\phi_{p',q}^0(\rad(C)),
\end{align*}
where $\rad(C)$ is radius of the open cluster $C$ at the origin,
that is, the maximum value of $n$ such that $0$ is joined by an open
path to the boundary of the box $[-n,n]^2$.
It follows that 
$$
\phi_{p',q}^0(\rad(C))<\oo,\qq p < \psd(q),
$$
is sufficient for $\pc(q)=\psd(q)$.

\begin{proof} 
We use a construction given in \cite{CC},
which was known earlier to one of the current
authors and to Paul Seymour. For odd $k$,
let $A_k$ be the event that $[0,2^k]\times[0,2^{k+1}]$ is
traversed long-ways. For even $k$, let $A_k$
be the event that $[0,2^{k+1}]\times[0,2^k]$ is
traversed long-ways.
By the positive-associativity and automorphism-invariance
of $\fpq^1$, 
under (\ref{finalcond}),
$$
\fpq^1\left(\bigcap_k A_k\right)\ge \prod_{k=1}^\oo \fpq^1(A_k)
> 0,\qq p>\psd(q).
$$
On the intersection of the $A_k$, there
exists an infinite open cluster, and therefore $\pc(q)\le \psd(q)$.
It is standard (see \cite{G03, G-RC})
that $\psd(q)\le \pc(q)$, and therefore equality holds
as claimed.
\end{proof}

\section{Acknowledgment} The work of the first author has
been
supported by the Engineering and Physical Sciences Research Council
through a PhD studentship.

\bibliography{references.bib}
\end{document}